**Floquet theory for integral and integro-differential equations.**

S. A. Belbas
Department of Mathematics
University of Alabama
Tuscaloosa, AL 35487-0350, USA.
*e-mail:  sbelbas@gmail.com*

Abstract:  We study the extension of Hill's method of infinite determinants to the case of integro-differential equations with periodic  coefficients and kernels. We develop the analytical theory  of such methods, and we obtain certain qualitative properties of the equations that determine the boundaries between  regions of dynamic stability and dynamic instability.

1. Introduction.

This paper contains the theoretical foundations (including analytical formulation of numerical methods) of extending the method of infinite determinants (Hill's determinants) to the case of integro-differential equations with periodic coefficients.

The standard Floquet theory, i.e. Floquet theory for systems of ordinary differential equations, has its origin in [F], and in the basic result about the existence of a monodromy matrix (see, e.g., [Av] for a succinct discussion of the concept of monodromy matrix) whose eigenvalues determine the stability or instability of the system; specifically, Schur-stable eigenvalues (i.e. eigenvalues having modulus less than 1) represent asymptotic stability, and those of modules greater than 1 correspond to instability. Thus, the standard Floquet theory involves sharp necessary and sufficient conditions for parametric stability, and an exact partition of the domain of system parameters into domains of asymptotic stability and instability. Some of the classical works on ODEs with periodic coefficients are [A, MW]. In this paper, we deal with the extension of Floquet theory to systems with memory. The theory of parametric stability and instability for integral and integro-differential equations is not a mere adaptation of the classical Floquet theory, but instead it involves new complications, raise new problems, and lead to new conditions, that have no counterpart in the theory of parametric stability and parametric resonance for ODEs.

In this paper, we are interested primarily in the extension of the method of infinite determinants to the case of integro-differential systems for the determination of the regions of parametric stability. The papers [ABDG, BBK, D2] have dealt with the existence of a monodromy matrix for periodic integro-differential systems; here, we derive in concise style the results about a monodromy matrix for the specific type of integro-differential systems (Volterra integro-differential systems over a time interval starting at $-\infty$). This sets up the scene for developing the method of infinite determinants. We should note that, in the theory of ODEs with periodic coefficients, there is a basic form of an ODE, due to Mathieu [M], that contains two parameters; many problems that arise in applications can be reduced to the Mathieu form, and consequently a plane diagram, on a plane with coordinates corresponding to those two parameters, suffices to depict the regions of dynamic stability and instability. In Floquet theory for ODE systems, a chart showing the regions of dynamic and instability, on a plane of two parameters of a scalar ODE, is known as a <u>Strutt chart</u>. We shall use the same terminology for the case of integro-differential system, although, of course, this version of "Strutt chart" is not the classical Strutt chart for the Mathieu-Hill equation. Low-order truncations of the determinantal equations give the parts of the Strutt chart that correspond to "small" values of the relevant parameters, and high-order truncations or other methods are needed for large values of the parameters. In the case of integro-differential equations, such a diagram is not possible, because the variety of Volterra kernels that could be utilized introduces large numbers of additional parameters. Consequently, the problems treated here fall into the general framework of multiparameter dynamic stability. In the ODE case, multiparameter problems have been treated in [SM].

The method of infinite determinants for ODE systems originated in the works of G. W. Hill on the motion of the Moon, specifically in [H1] and in a booklet privately published by G. W. Hill; see also [W] and [H2]. Hill treated infinite linear systems as if they were finite systems, and utilized conditions involving infinite determinants, for which he then took finite truncations in order to perform actual calculations. Hill's non-rigorous work gained the admiration of Poincaré, and it was later made rigorous primarily through the investigations of Poincaré and H. von Koch. Research on infinite determinants has continued unto the present day (see, e.g., [GLMZ, JZ]). A general account of current problems in dynamic stability and parametric resonance for systems governed by differential equations is provided by [FN].

Our interest in Floquet theory for integro-differential systems is motivated from problems of dynamic stability for viscoelastic systems with constitutive laws involving general Volterra integral operators (as distinct from the differential or the complex models of viscoelastic constitutive laws, which can be expressed via differential operators). This in turn is motivated by the ubiquitous appearance of viscoelastic materials in many applications, for instance, food materials, biological materials, synthetic materials for biomedical and robotic applications, structural elements including elements for the damping of vibrations. For a general review of classical viscoelastic models, see [C]; for models involving general Volterra integral operators, see [D]. The method of infinite determinants has been used extensively for the dynamic stability of elastic systems, for example in [Bo1, Bo2]. In this paper, we extend the method of infinite determinants to integro-differential applications. This extension creates new problems, new techniques, and new results about the boundaries of the regions of dynamic stability and instability in the space of suitable parameters.

As stated in the beginning of this section, inn this paper we study the analytical theory of infinite determinants for integro-differential Floquet theory. We shall present explicit problem of viscoelastic vibrations, and actual calculations, in future papers.

## 2. General properties of periodic integro-differential systems.

We are interested in developing a theory of periodic solutions for periodic integro-differential systems that will extend the classical Floquet theory of periodic differential systems. The appropriate form of integro-differential systems that admits a Floquet-type theory turns out to be

$$\frac{dx(t)}{dt} = A(t)x(t) + \int_{-\infty}^{t} K(t,s)x(s)\,ds \tag{2.1}$$

where $x(t)$ takes values in $\mathbb{R}^n$, $A(t)$ is an $n \times n$ matrix-valued function that is continuous for $t \in \mathbb{R}$ and periodic with period T, i.e. $A(t+T) = A(t)\ \forall t \in \mathbb{R}$, and $K(t, s)$ is continuous for $(t,s) \in D := \{(\tau,\sigma) : -\infty < \sigma \leq \tau < +\infty\}$ and satisfies $K(t+T, s+T) = K(t,s)\ \forall (t,s) \in D$ and $\iint_D |K(t,s)|\,dt\,ds < \infty$. We introduce the name *simultaneous periodicity* for the property $K(t+T, s+T) = K(t,s)$. For example, every kernel of convolution type, i.e. $K(t,s) = \varphi(t-s)$, has the property of simultaneous periodicity, although it need not be separately periodic in either t or s. Many other types of more general kernels also have the property of simultaneous periodicity. The question of existence and uniqueness of solutions of (2.1) for every initial condition $x(t_0) = x_0$ is studied in Appendix 1 of the present paper; for the purposes of the current section, we assume existence and uniqueness of the Cauchy problem for (2.1) for every initial data $(t_0, x_0)$. Our first task is to demonstrate the Floquet-like properties of (2.1).

First, by the linearity of (2.1), every solution of (2.1) can be expressed as a linear combination of solutions $x_1(t), x_2(t), ..., x_n(t)$ where each $x_i(t)$ solves (2.1) with initial condition $x_i(t_0) = x_{0,i}$ and the vectors $\{x_{0,i} : 1 \leq i \leq n\}$ are linearly independent (for example,

$x_{0,i} = e_i = (i^{\text{th}}$ column of the $n \times n$ unit matrix) ).

Second, if $x(t)$ is a solution of (2.1), then $x(t+T)$ is also a solution of the same equation. Indeed, the only term that is different from the classical Floquet theory is the integral term, and for that term we have, by dint of the simultaneous periodicity property of K,

$$\int_{-\infty}^{t+T} K(t+T,s)x(s)\,ds \underset{(\text{set } s'=s-T)}{=} \int_{-\infty}^{t} K(t+T,s'+T)x(s'+T)\,ds' = \\ = \int_{-\infty}^{t} K(t,s')x(s'+T)\,ds' \tag{2.2}$$

Then, as in the classical case of Floquet theory, we conclude that there exists a monodromy matrix, say M, such that

$$x(t+T) = Mx(t). \tag{2.3}$$

The stability or instability of the solutions of (2.1) depends on the spectral properties of M. As in the standard Floquet theory for ODEs, the eigenvalues of M and the characteristic and minimal polynomials of M do not change when we consider another solution, say y(t), of (2.1), although the matrix M corresponding to x(t) is not necessarily the same with the monodromy matrix for y(t). Consequently, if $X(t,t_0)$ satisfies

$$\frac{\partial X(t,t_0)}{\partial t} = A(t)X(t,t_0) + \int_{-\infty}^{t} K(t,s)X(s,t_0)ds \; ; \; X(t_0,t_0) = I_n \, (= \text{the } n \times n \text{ unit matrix}) \tag{2.4}$$

then it suffices to consider the matrix $M_X := X(t_0 + T, t_0)$. The solutions of (2.1) will be stable if every (generally complex) eigenvalue $\lambda$ of $M_X$ satisfies $|\lambda| \leq 1$ and all eigenvalues with $|\lambda| = 1$ are simple eigenvalues, and unstable otherwise.

These conclusions require that $\int_{-\infty}^{t} K(t,s)\exp(\lambda(t-s))ds$ be absolutely convergent for all real values of $\lambda$. This will be the case, for example, if the kernel K(t, s) satisfies

$$|K(t,s)| \leq C_K \exp(-\mu(t-s)^\beta), \; \beta > 1. \tag{2.5}$$

For a simple (generally complex) non-zero eigenvalue $\lambda$, let $\Lambda$ be one branch of (the complex) $\ln \lambda$. If $x_1(t)$ is a solution that satisfies $x_1(t+T) = \lambda x_1(t)$, then $y_1(t) := \exp\left(-\frac{\Lambda t}{T}\right)x_1(t)$ satisfies $y_1(t+T) = y_1(t)$, and therefore $x_1$ has the interpretation $x_1(t) = \exp\left(\frac{\Lambda t}{T}\right)y_1(t)$ where $y_1$ is periodic in t with period T. An analogous conclusion holds if a solution $x_2(t)$ satisfies $x_2(t+T) = -\lambda x_2(t)$: in that case, the function $y_2(t) := \exp\left(-\frac{\Lambda t}{T}\right)x_2(t)$ is antiperiodic with antiperiod T, i.e. it satisfies $y_2(t+T) = -y_2(t)$, and $x_2(t)$ has the representation $x_2(t) = \exp\left(\frac{\Lambda t}{T}\right)y_2(t)$ with antiperiodic $y_2$.

Also, if some solution $x_3(t)$ satisfies $x_3(t+T) = Mx_3(t)$ and the $n \times n$ matrix M has the representation $M = \exp(\mathfrak{m})$, then $x_3$ can be represented as $x_3(t) = \exp\left(\frac{\mathfrak{m}t}{T}\right)y_3(t)$ where $y_3$ is perioduc with period T, whereas if another solution $x_4(t)$ satisfies $x_4(t+T) = -Mx_4(t)$, then $x_4$ is representable as $x_4(t) = \exp\left(\frac{\mathfrak{m}t}{T}\right)y_4(t)$ where $y_4$ is antiperiodic with antiperiod T.

Of course, antiperiodic solutions can also be handled by taking the (complex) logarithm of $\lambda$ or $\mathfrak{m}$, i.e. by adding $i\pi$ to $\Lambda$ or $i\pi I_n$ to M; it is nevertheless desirable to also have an explicit representation of solutions of the integrodifferential system in terms of antiperiodic functions.

These claims are easily proved in essentially the same way as in the standard Floquet theory for ODEs, and we do not dwell on the proofs here.

## 3. A second order integrodifferential system.

Because of its relevance for applications, we consider the following system with one-dimensional unknown function x(t):

$$\ddot{x}(t) + a(t)x(t) = \int_{-\infty}^{t} K(t,s)x(s)\,ds \tag{3.1}$$

where the function a(t) is continuous and periodic with period T, and the one-dimensional kernel K is continuous, jointly periodic in (t, s) with period T, and satisfies (3.1).

We take a(t) in the form

$$a(t) = a_0 + a_1 \cos(\vartheta t)\,;\quad \vartheta = \frac{2\pi}{T} \tag{3.2}$$

In this case, the monodromy matrix is a $2\times 2$ matrix and thus its characteristic equation is a quadratic equation, say $\lambda^2 + p\lambda + q = 0$ with real coefficients p and q. The necessary and sufficient condition for both roots of the characteristic equation to satisfy $|\lambda|<1$ is $\{|q|<1 \text{ and } |p|<1+q\}$. For the purpose of characterizing the boundaries between dynamic stability and dynamic instability, the case of one root that satsfies $|\lambda|=1$ does not need to be considered separately. On the pq-plane, the domain of all (p, q) that satisfy $\{|q|<1 \text{ and } |p|<1+q\}$ is the interior of a triangle with vertices A, B, C, where
$A = (p_A, q_A) = (0, -1)$, $B = (p_B, q_B) = (2, 1)$, $C = (p_C, q_C) = (-2, 1)$.

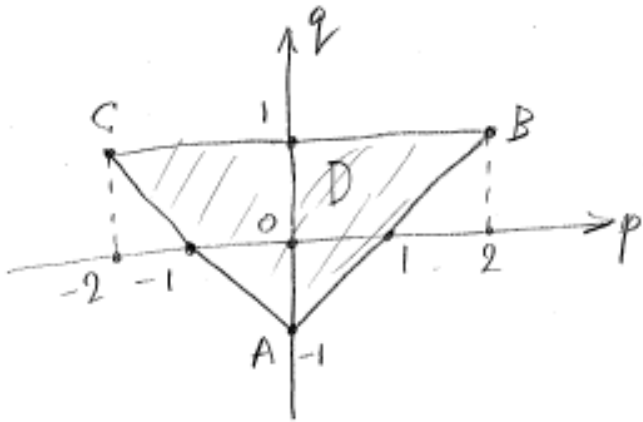

*Figure 2.1.*
Stability domain on the pq-plane.
--------------------------------------------------------------------------------

Thus the boundary of the domain of dynamic stability on the pq-plane consists of the 3 open straight-line segments

$]AB[ := \{(p,q) : p = 1+q, -1 < q < 1\}$, $]AC[ := \{(p,q) : p = -(1+q), -1 < q < 1\}$,
$]BC[ := \{(p,q) : q = 1, -2 < p < 2\}$,

and the 3 vertices A, B, C.

For solutions of (2.1) that correspond to values of the parameters $a_0, a_1$, and $\vartheta$ on the boundary of the domains of stability (i.e. values of $a_0, a_1$, and $\vartheta$ for which the coefficients of the characteristic equation of the monodromy matrix belong to the boundary of the stability domain on the pq-plane), we have six possibilities. The first two possibilities concern the existence of solutions x(t) that have the form

$$x(t) = \exp(i\omega' t) y(t); \quad 0 \leq \omega \leq \pi; \omega' := \frac{\omega}{T} \quad , \tag{3.3}$$

where y(t) is either periodic with period T, in which case we take y(t) in the form

$$y(t) = \sum_{n=-\infty}^{\infty} Y_n \exp(in\vartheta t), \tag{3.4}$$

or antiperiodic with antiperiod T, in the form

$$y(t) = \sum_{n=-\infty}^{\infty} Z_n \exp\left(\tfrac{1}{2} i(2n+1)\vartheta t\right). \tag{3.5}$$

We state here that, in the context of the present paper, all Fourier and antiperiodic trigonometric series are intended to represent real-valued functions, thus, in expressions like (3.4) and (3.5), it is automatically postulated that $Y_{-n} = \overline{Y}_n$, $Z_{-(n+1)} = \overline{Z}_n$, respectively, for all integers n, where the overbar signifies the operation of taking the conjugate of a complex number. All Fourier series of the types (3.4) and (3.5) can also be written as Fourier or antiperiodic trigonometric series with sine and cosine terms and with real coefficients, thus, for instance,

$$\sum_{n=-\infty}^{\infty} Y_n \exp(in\vartheta t) = A_0 + \sum_{n=1}^{\infty} (A_n \cos(n\vartheta t) + B_n \sin(n\vartheta t));$$
$$A_0 = Y_0; A_n = Y_n + Y_{-n}, B_n = i(Y_n - Y_{-n}), \text{ if } n \neq 0 .$$

and similarly for the anti-periodic Fourier-like series in (3.5).

The remaining possibilities correspond to the simultaneous existence of two solutions $x_1(t)$, $x_2(t)$, of which $x_1(t)$ is either periodic with period T or antiperiodic with antiperiod T, and $x_2(t) = e^{-\lambda t} y(t)$ ($\lambda > 0$) where y(t) is either periodic of the form (3.4) or antiperiodic of the form (3.5).

It must be emphasized that the occurrence of all these possibilities is peculiar to integro-differential systems, whereas only the first two possibilities (periodic or anti-periodic solutions) with either $\omega = 0$ or $\omega = \pi$, are relevant to the corresponding ODE systems. This is due to the fact that, in the absence of the integral terms, i.e. in the case of a standard Mathieu-Hill second order ODE, it is shown in many treatises (including, for instance, [Bo]) that the constant q in the characteristic equation can be taken equal to 1, and then a continuity argument shows that the boundaries between regions of stability and regions on instability (in the space of the parameters of the ODE) are obtained by considering periodic and antiperiodic solutions of the Mathieu-Hill ODE. With q=1, we are limited to the line BC in Fig. 2.1, and the boundaries between domains of stability and domains of instability correspond to complex roots of the characteristic equation and thus to periodic and antiperiodic solutions. These considerations and arguments do not carry over to the case of integro-differential equations, and thus we are obliged to consider more complicated possibilities, i.e. all 3 sides of the triangle ABC of Fig. 2.1.

Among many different possible ways of describing the boundaries between stability and instability, we shall have to make different choices in different cases. This will become clear as we examine each case separately. The traditional way (for the Mathieu-Hill ODE) to describe boundaries between regions of dynamic stability and dynamic instability is to draw diagrams (the so-called Strutt - charts) on the plane of the parameters $a_1$ and $\vartheta$ (or something similar, like, e.g., $\frac{a_1}{a_0}$ and $\vartheta$); however, for the questions that arise for integro-differential equations, the traditional choices are not always the best.

We start by examining the first two possibilities. For a solution of the form (3.3) with y satisfying (3.4), we set

$$K_{nm}(\vartheta, \gamma) := \tfrac{1}{T} \int_{t=0}^{T} \int_{s=-\infty}^{t} K(t,s) \exp(im\vartheta s - in\vartheta t + \gamma(s-t)) \, ds \, dt . \tag{3.6}$$

Then the existence of a solution of the form (3.3) with y expressed as in (3.4) implies

$$-(n\vartheta + \omega')^2 Y_n + a_0 Y_n + \tfrac{1}{2} a_1 (Y_{n-1} + Y_{n+1}) = \sum_{m=-\infty}^{\infty} K_{nm}(\vartheta, i\omega') Y_m . \tag{3.7}$$

Similarly, for a solution of the form (3.3) with antiperiodic y(t) given by (3.5), we set

$$\widetilde{K}_{nm}(\vartheta, \gamma) := \tfrac{1}{T} \int_{t=0}^{T} \int_{s=-\infty}^{t} K(t,s) \exp\!\left(\tfrac{i}{2}((2m+1)\vartheta s - (2n+1)\vartheta t) + \gamma(s-t)\right) ds$$

(3.8)

and we find that the coefficients of the antiperiodic function in (3.5) must satisfy

$$-\left(\omega'+\tfrac{2n+1}{2}\vartheta\right)^2 Z_n + a_0 Z_n + \tfrac{1}{2}a_1(Z_{n-1}+Z_{n+1}) = \sum_{m=-\infty}^{\infty} \widetilde{K}_{nm}(\vartheta, i\omega') Z_m \qquad (3.9)$$

We define the doubly infinite matrices $\mathbf{D}(\vartheta, \gamma)$ and $\widetilde{\mathbf{D}}(\vartheta, \gamma)$ by

$$\begin{aligned}
&\mathbf{D}(\vartheta,\gamma) = [D_{nm}(\vartheta,\gamma)]_{(n,m)=(-\infty,-\infty)}^{(\infty,\infty)} ; \\
&D_{nm}(\vartheta,\gamma) = \delta_{nm} - \tfrac{1}{(in\vartheta+\gamma)^2}\left\{K_{nm}(\vartheta,\gamma) - a_0\delta_{nm} - \tfrac{1}{2}a_1\delta_{n-1,m} - \tfrac{1}{2}a_1\delta_{n+1,m}\right\}; \\
&\widetilde{\mathbf{D}}(\vartheta,\gamma) = [\widetilde{D}_{nm}(\vartheta,\gamma)]_{(n,m)=(-\infty,-\infty)}^{(\infty,\infty)} ; \\
&\widetilde{D}_{nm}(\vartheta,\gamma) = \delta_{nm} - \tfrac{1}{(\tfrac{2n+1}{2}i\vartheta+\gamma)^2}\left\{\widetilde{K}_{nm}(\vartheta,\gamma) - a_0\delta_{nm} - \tfrac{1}{2}a_1\delta_{n-1,m} - \tfrac{1}{2}a_1\delta_{n+1,m}\right\}
\end{aligned} \qquad (3.10)$$

where $\delta_{ij}$ stands for Kronecker's delta. The definition of $D_{0,m}(\vartheta,0)$ is an exception. In that case, the row $[D_{0,m}(\vartheta,0)]$ is replaced by

$$\hat{D}_{0,m}(\vartheta,0) := -K_{0,m}(\vartheta,0) + a_0\delta_{0,m} + \tfrac{1}{2}a_1(\delta_{-1,m} + \delta_{1,m}) . \qquad (3.11)$$

It will be shown a little further down that, under certain conditions on the kernel K, the matrices $\mathbf{D}$ and $\widetilde{\mathbf{D}}$ satisfy one of the sufficient conditions that ensure that the existence of non-zero solutions of the systems (3.7) and (3.9) is equivalent to the vanishing of the determinants det $(\mathbf{D}(\vartheta, i\omega))$ and det $(\widetilde{\mathbf{D}}(\vartheta, i\omega))$, respectively, and that those infinite determinants converge and can be approximated by finite determinants in the sense that the value of each infinite determinant is obtained as the limit of a sequence of finite determinants:

$$\det(\mathbf{D}(\vartheta, i\omega')) = 0, \qquad (3.12)$$

$$\det(\widetilde{\mathbf{D}}(\vartheta, i\omega')) = 0. \qquad (3.13)$$

Along the side ]AB[ , in Fig. 2.1, roots of the characteristic equation occur as conjugate pairs, and consequently the conditions that correspond to [AB[ are

$$\det(\mathbf{D}(\vartheta, i\omega')) = \det(\mathbf{D}(\vartheta, -i\omega')) = 0; \quad \det(\widetilde{\mathbf{D}}(\vartheta, i\omega')) = \det(\widetilde{\mathbf{D}}(\vartheta, -i\omega')) = 0 . \qquad (3.14)$$

We note that the quantities $K_{nm}(\vartheta, i\omega)$ have the property $\overline{K}_{nm}(\vartheta, i\omega) = K_{-n,-m}(\vartheta, -i\omega)$ (the overbar denotes complex conjugation), and consequently the matrix $\mathbf{D}(\vartheta, -i\omega')$ is the Hermitian-adjoint (transposition of the conjugate) of $\mathbf{D}(\vartheta,i\omega')$. Thus the first system in (3.14) is tantamount to

$$\det(\mathbf{D}(\vartheta,i\omega')) = 0. \tag{3.15}$$

Similarly, the second system in (3.14) can be replaced by

$$\det(\widetilde{\mathbf{D}}(\vartheta,i\omega')) = 0. \tag{3.16}$$

The remaining cases (sides AB and AC in Fig. 2.1) lead to the systems

$$\begin{aligned}&\det(\mathbf{D}(\vartheta,0)) = \det(\mathbf{D}(\vartheta,-\lambda)) = 0 \,;\; \det(\widetilde{\mathbf{D}}(\vartheta,0)) = \det(\widetilde{\mathbf{D}}(\vartheta,-\lambda)) = 0;\\ &\det(\mathbf{D}(\vartheta,0)) = \det(\widetilde{\mathbf{D}}(\vartheta,-\lambda)) = 0 \,;\; \det(\widetilde{\mathbf{D}}(\vartheta,0)) = \det(\mathbf{D}(\vartheta,-\lambda)) = 0\end{aligned}. \tag{3.17}$$

The two systems in (3.15, 3.16) can be put in real form. We denote by $K_{nm,C}(\vartheta,\omega)$ and $K_{nm,S}(\vartheta,\omega)$, respectively, the real and imaginary parts of $K_{nm}(\vartheta,i\omega)$. In other words,

$$K_{nm,C}(\vartheta,\omega) := \frac{1}{T}\int_{t=0}^{T}\int_{s=-\infty}^{t} K(t,s)\cos(\vartheta(ms-nt)+\omega(s-t))\,ds\,dt;$$

$$K_{nm,S}(\vartheta,\omega) := \frac{1}{T}\int_{t=0}^{T}\int_{s=-\infty}^{t} K(t,s)\sin(\vartheta(ms-nt)+\omega(s-t))\,ds\,dt\,.$$

Analogous definitions hold for $\widetilde{K}_{nm}(\vartheta,i\omega)$.

The matrix $\mathbf{D}(\vartheta,i\omega)$ can be expressed as $\mathbf{D}(\vartheta,i\omega) = \mathbf{D}_C(\vartheta,\omega) + i\mathbf{D}_S(\vartheta,\omega)$, where $\mathbf{D}_C(\vartheta,\omega)$ and $\mathbf{D}_S(\vartheta,\omega)$ are obtained from $\mathbf{D}(\vartheta,i\omega)$ by replacing $K_{nm}(\vartheta,i\omega)$ by $K_{nm,C}(\vartheta,\omega)$ and $K_{nm,S}(\vartheta,\omega)$, respectively.

We shall use the superscript (N) to denote the truncation of an infinite determinant to an N-dimensional determinant. In order to utilize the well-known formulae for the determinant of a sum of two finite dimensional $N \times N$ matrices, we consider all N-dimensional multi-indices $\alpha$ with entries 0 or 1 (thus there are $2^N$ such multi-indices), and we define $\mathbf{D}^{(N)(\alpha)}(\vartheta,\omega')$ to be the determinant that has k-th row equal to the k-th row of $\mathbf{D}_C(\vartheta,\omega)$ if the k-th component of $\alpha$ is 0, and equal to th k-th row of $\mathbf{D}_S(\vartheta,\omega)$ if the k-th component of $\alpha$ is 1. We use the notation $|\alpha|$ for the sum of all components of the multi-index $\alpha$. Then, by the well known theorems about the determinant of the sum of two matrices, plus the standard properties of determinants, we have

$$\det(\mathbf{D}^{(N)}(\vartheta, i\omega')) = \sum_{\alpha:\,|\alpha|\equiv 0\bmod 4} \det(\mathbf{D}^{(N)(\alpha)}(\vartheta,\omega')) - \sum_{\alpha:\,|\alpha|\equiv 2\bmod 4} \det(\mathbf{D}^{(N)(\alpha)}(\vartheta,\omega')) +$$
$$+ i\left[\sum_{\alpha:\,|\alpha|\equiv 1\bmod 4} \det(\mathbf{D}^{(N)(\alpha)}(\vartheta,\omega')) - \sum_{\alpha:\,|\alpha|\equiv 3\bmod 4} \det(\mathbf{D}^{(N)(\alpha)}(\vartheta,\omega'))\right] \quad .$$

(3.18)

Consequently, the N-dimensional truncation of the system with complex coefficients (3.15) is equivalent to the system with real coefficients

$$\sum_{\alpha:\,|\alpha|\equiv 0\bmod 4} \det(\mathbf{D}^{(N)(\alpha)}(\vartheta,\omega')) - \sum_{\alpha:\,|\alpha|\equiv 2\bmod 4} \det(\mathbf{D}^{(N)(\alpha)}(\vartheta,\omega')) = 0;$$
$$\sum_{\alpha:\,|\alpha|\equiv 1\bmod 4} \det(\mathbf{D}^{(N)(\alpha)}(\vartheta,\omega')) - \sum_{\alpha:\,|\alpha|\equiv 3\bmod 4} \det(\mathbf{D}^{(N)(\alpha)}(\vartheta,\omega')) = 0 \quad .$$

(3.19)

In (3.15), the domain of $\lambda$ is $[0,+\infty)$. For numerical implementations, we may take $\lambda = \tan\varphi$, $\varphi \in [0,\frac{\pi}{2})$, and carry out the calculations with a uniform-mesh discretization of $[0,\frac{\pi}{2})$.

We note here that, in the literature on parametric resonance for ODEs, the <u>normality</u> condition (due to H. von Koch) for infinite determinants is most well known; however, there exist also more general sufficient conditions for infinite matrices and their determinants to have the required properties (as well as <u>necessary and</u> sufficient conditions, which, however, we shall not use here). For an infinite matrix with elements $M_{nm} = \delta_{nm} + S_{nm}$, <u>normality</u> is the condition $\sum_{n,m} |S_{nm}| < \infty$. A more general sufficient condition is $\sum_{n} \sup_{m} |S_{nm}| < \infty$; see [K], and section 6.5 of [P]. Under the condition that the kernel K is bounded over the domain $\{s \le 0,\, 0 \le t \le T\} \cup \{0 \le s \le t \le T\}$, the matrices $\mathbf{D}(\vartheta, i\omega)$ and $\widetilde{\mathbf{D}}(\vartheta, i\omega)$ satisfy this latter sufficient condition, because coefficients $\frac{1}{(in\vartheta+i\omega)^2}$ and $\frac{1}{(\frac{2n+1}{2}i\vartheta+i\omega)^2}$ are terms of absolutely convergent series as n runs through the integers.

For the remaining matrices, shown in (3.15), the condition of faster-than-exponential rate of decay of the kernel K, shown in Eq. (2.5), together with form of the matrices defined in (3.10), again ensures that the second sufficient condition mentioned above is satisfied for all real values of $\gamma$; for non-positive values of $\gamma$, the condition of exponential order of decay suffices.

## 4. The numerical evaluation of the terms associated with the memory kernel.

The terms

$$K_{nm}(\vartheta,\gamma) = \tfrac{1}{T}\int_{t=0}^{T}\int_{s=-\infty}^{t} K(t,s)\exp(i\vartheta(ms-nt)+\gamma(s-t))\,ds\,dt$$

may be evaluated numerically by a combination of Gauss-Legendre approximate integration and an approximate integration formula over a finite interval (for example a Newton-Cotes method).

By using variables t and $\xi := t-s$, we transform the above integral to

$$K_{nm}(\vartheta,\gamma) = \tfrac{1}{T}\int_{t=0}^{T}\int_{\xi=0}^{\infty} K(t,t-\xi)\exp(i\vartheta(m-n)t-(i\vartheta m+\gamma)\xi)\,d\xi\,dt\,.$$

Now, an integral of the form $\int_0^\infty e^{-x} f(x)\,dx$ can be approximated via Gauss-Legendre approximation by $\sum_{i=1}^{N} w_i f(x_i)$ where $w_i$ and $x_i$ are, respectively, the weights and nodes of Gauss-Legendre quadrature. $\int_0^\infty e^{-vx} f(x)\,dx$ is approximated by $\sum_{i=1}^{N} \frac{w_i}{v} f\left(\frac{x_i}{v}\right)$. Thus an integral $\int_0^\infty g(x)\,dx$, where g satisfies $|g(x)| \le C_g \exp(-vx)\ \forall x>0$, is approximated by $\sum_{i=1}^{N} e^{vx_i}\left(\frac{w_i}{v}\right)g\left(\frac{x_i}{v}\right)$. On the other hand, an integral $\int_0^T h(y)\,dy$ can be approximated, via a Newton-Cotes formula, by an expression of the form $\sum_{j=1}^{M} u_j h(y_j)$. Consequently, $K_{nm}(\vartheta,\beta)$ can be approximated by

$$\sum_{i=1}^{N}\sum_{j=1}^{M} u_j \exp(vx_i)\left(\tfrac{w_i}{v}\right) K\!\left(y_j, y_j - \tfrac{x_i}{v}\right)\exp\!\left(i\vartheta(m-n)y_j - \tfrac{i\vartheta m+\gamma}{v}x_i\right).$$

The value of $v$ is taken as $v = \mu + \Re\gamma$, and the approximation is valid for $v > 0$.

Similarly, $\widetilde{K}_{nm}(\vartheta,\gamma)$ can be approximated by

$$\sum_{i=1}^{N} \sum_{j=1}^{M} u_j \exp(vx_i) \left(\frac{w_i}{v}\right) K\left(y_j, y_j - \frac{x_i}{v}\right) \exp\left(i\vartheta(m'-n')y_j - \frac{i\vartheta m'+\gamma}{v} x_i\right);$$

$$m' := m + \tfrac{1}{2}, \quad n' := n + \tfrac{1}{2}$$

These formulae provide all the necessary analytical information for the approximate evaluation of the integrals appearing in the relevant equations for the determination of stable values of the parameters of the problem of section 3.

5. The truncated determinantal systems for the case of periodic solutions.

The methods of obtaining concrete approximations to the boundaries of the regions of dynamic stability require treatment of systems obtained by truncating the infinite determinants of the previous section. As in the case of the analogous problems for ordinary differential equations, the first results, that are tractable by analytical approximations, are obtained by using 3-dimensional approximations to the infinite determinants.

For a $1 \times 1$ approximation of the determinantal equations for the first case, viz. (3.11), we find the equations

$$a_0 - K_{00}(\vartheta, i\omega') - \omega'^2 = 0 \tag{5.1}$$

which gives the range of the quasi-static resonance values

$$a_0^*(\vartheta, \omega') = \omega'^2 + K_{00}(\vartheta, i\omega'), \quad 0 \leq \omega \leq \tfrac{\pi}{2} \tag{5.2}$$

Only the values of $\omega'$ that give real values for the right-hand side of (5.2) are physically meaningful.

We shall call the values given by (5.2) <u>quasi-static</u> values because they correspond to values of only the time-independent part of a(t), i.e. values of $a_0$, and they do not involve the values of $a_1$ which is the oscillatory part of a(t).

We note that $K_{00}(\vartheta, i\omega')$ is given by

$$K_{00}(\vartheta, i\omega') = \tfrac{1}{T}\int_{t=0}^{T}\int_{s=-\infty}^{t} K(t,s)\exp(i\omega'(s-t))\,ds\,dt =$$
$$= \tfrac{1}{T}\int_{t=0}^{T}\int_{s=-\infty}^{t} K(t,s)\cos(\omega'(s-t))\,ds\,dt + \tfrac{i}{T}\int_{t=0}^{T}\int_{s=-\infty}^{t} K(t,s)\sin(\omega'(s-t))\,ds\,dt \tag{5.3}$$

and consequently the values of $\omega'$ that give real values in (4.2) are solutions of the equation

$$\int_{t=0}^{T}\int_{s=-\infty}^{t} K(t,s)\sin(\omega'(s-t))\,ds\,dt = 0 . \tag{5.4}$$

The expression in (5.4) can be transformed via the substitution $\xi = t - s$, $\eta = t$, which was used in the previous section, so that

$$\int_{t=0}^{T} \int_{s=-\infty}^{t} K(t,s)\sin(\omega'(s-t))\,ds\,dt = -\int_{\eta=0}^{T} \int_{\xi=0}^{+\infty} K(\eta,\eta-\xi)\sin(\omega'\xi)\,d\xi\,d\eta \ .$$

By setting

$$K^{\wedge}(\xi) := \int_{\eta=0}^{T} K(\eta,\eta-\xi)\,d\eta$$

we can write (5.4) in the form

$$\int_{\xi=0}^{+\infty} K^{\wedge}(\xi)\sin(\omega'\xi)\,d\xi = 0 \qquad (5.5)$$

Eq. (5.5) is an equation about the zeroes of Fourier integrals. Eq. (5.5), however, forms a different problem from that of the general question of the distribution of zeros of Fourier integrals: in (5.5), we are searching only for real roots, when the variable ω is restricted to the interval $[0,\frac{\pi}{2}]$ . A first approximation to the real roots may be obtained by using a series expansion for the left-hand side of (5.5), namely

$$\int_{\xi=0}^{+\infty} K^{\wedge}(\xi)\sin(\omega'\xi)\,d\xi = \sum_{k=0}^{\infty} \frac{(-1)^k \omega'^{2k+1}}{(2k+1)!} \int_{\xi=0}^{+\infty} K^{\wedge}(\xi)\xi^{2k+1}\,d\xi . \qquad (5.6)$$

By keeping the first 3 terms of the series expansion above, we obtain the following approximation to (5.5):

$$A_0 - A_1\omega'^2 + A_2\omega'^4 = 0 \ ; \quad A_k := \tfrac{1}{(2k+1)!}\int_{\xi=0}^{+\infty} K^{\wedge}(\xi)\xi^{2k+1}\,d\xi \qquad (5.7)$$

with solutions (provided the expression below gives real values for ω')

$$\omega'^* = \sqrt{\frac{A_1 \pm \sqrt{A_1^2 - 4A_0 A_2}}{2A_2}} \qquad (5.8)$$

We shall call the value(s) $\omega'^*$ above <u>the quasi-static frequency(ies) of parametric resonance</u>.

The $3\times 3$ approximation to the infinite determinantal equation (3.11) is

$$\begin{vmatrix} K_{-1,-1}(\vartheta,i\omega') + (\omega'-\vartheta)^2 - a_0 & K_{-1,0}(\vartheta,i\omega') - \tfrac{1}{2}a_1 & K_{-1,1}(\vartheta,i\omega') \\ K_{0,-1}(\vartheta,i\omega') - \tfrac{1}{2}a_1 & K_{00}(\vartheta,i\omega') + \omega'^2 - a_0 & K_{01}(\vartheta,i\omega') - \tfrac{1}{2}a_1 \\ K_{1,-1}(\vartheta,i\omega') & K_{10}(\vartheta,i\omega') - \tfrac{1}{2}a_1 & K_{11}(\vartheta,i\omega') + (\omega'+\vartheta)^2 - a_0 \end{vmatrix} = 0$$

(5.9)

This is a quadratic equation for $a_1$, viz.

$$C_0 + C_1 a_1 + C_2 a_1^2 = 0 \tag{5.10}$$

where

$$\begin{aligned}
C_0 &= -\left(K_{00}(\vartheta,i\omega') + \omega'^2 - a_0\right)\left(K_{-1,-1}(\vartheta,i\omega') + (\omega'-\vartheta)^2 - a_0\right)\left(K_{11}(\vartheta,i\omega') + (\omega'+\vartheta)^2 - a_0\right) - \\
&\quad - K_{-1,0}(\vartheta,i\omega')K_{01}(\vartheta,i\omega')K_{1,-1}(\vartheta,i\omega') - K_{0,-1}(\vartheta,i\omega')K_{1,0}(\vartheta,i\omega')K_{-1,1}(\vartheta,i\omega') + \\
&\quad + K_{-1,0}(\vartheta,i\omega')K_{0,-1}(\vartheta,i\omega')\left(K_{11}(\vartheta,i\omega') + (\omega'+\vartheta)^2 - a_0\right) + \\
&\quad + K_{01}(\vartheta,i\omega')K_{10}(\vartheta,i\omega')\left(K_{-1,-1}(\vartheta,i\omega') + (\omega'-\vartheta)^2\right); \\
C_1 &= \tfrac{1}{2}\Big[ K_{-1,0}(\vartheta,i\omega')K_{1,-1}(\vartheta,i\omega') + K_{01}(\vartheta,i\omega')K_{1,-1}(\vartheta,i\omega') + K_{01}(\vartheta,i\omega')K_{-1,1}(\vartheta,i\omega') + \\
&\quad + K_{0,-1}(\vartheta,i\omega')K_{-1,1}(\vartheta,i\omega') - \left(K_{01}(\vartheta,i\omega') + K_{10}(\vartheta,i\omega')\right)\left(K_{-1,-1}(\vartheta,i\omega') + (\omega'-\vartheta)^2 - a_0\right) - \\
&\quad - \left(K_{-1,0}(\vartheta,i\omega') + K_{0,-1}(\vartheta,i\omega')\right)\left(K_{11}(\vartheta,i\omega') + (\omega'+\vartheta)^2 - a_0\right) \Big]; \\
C_2 &= \tfrac{1}{4}\Big[ K_{11}(\vartheta,i\omega') + K_{-1,-1}(\vartheta,i\omega') - K_{1,-1}(\vartheta,i\omega') - K_{-1,1}(\vartheta,i\omega') \Big] + \tfrac{1}{2}\left(\omega'^2 + \vartheta^2 - a_0\right)
\end{aligned}$$

(5.11)

Our exposition below will be facilitated by introducing the following notation. With the convention that $\Re z$, $\Im z$ will stand for the real and imaginary parts, respectively, of a complex number z, we set

$$Q_1 := \begin{vmatrix} \Re C_2 & \Re C_0 \\ \Im C_2 & \Im C_0 \end{vmatrix}^2 - \begin{vmatrix} \Re C_2 & \Re C_1 \\ \Im C_2 & \Im C_1 \end{vmatrix} \begin{vmatrix} \Re C_1 & \Re C_0 \\ \Im C_1 & \Im C_0 \end{vmatrix}, \quad Q_2 := (\Re C_1)^2 - 4(\Re C_0)(\Re C_2).$$

(5.12)

The necessary and sufficient condition for (5.10) to have at least one real root, provided that $\Re C_2 \neq 0$ and $\Im C_2 \neq 0$, is that $Q_2 \geq 0$ and $Q_1 = 0$. This follows from standard elimination theory.

Further information about real roots of polynomials with complex coefficients may be found in [Ba]. There are known conditions for a polynomial with complex coefficients to have exactly a specified number of real roots, but we shall not repeat those conditions here.

## 6. Truncated determinants for the case of antiperiodic solutions.

Below, we write simply $\widetilde{K}_{nm}$ for $\widetilde{K}_{nm}(\vartheta, i\omega')$.

For the case of antiperiodic solutions, the first meaningful approximation of the infinite-determinantal equation (3.12) is an equation with a $2 \times 2$ determinant:

$$\begin{vmatrix} a_0 - (\omega' + \tfrac{1}{2}\vartheta)^2 - \widetilde{K}_{-1,-1} & \tfrac{1}{2}a_1 - \widetilde{K}_{-1,0} \\ \tfrac{1}{2}a_1 - \widetilde{K}_{0,-1} & a_0 - (\omega' + \tfrac{1}{2}\vartheta)^2 - \widetilde{K}_{00} \end{vmatrix} = 0 . \tag{6.1}$$

In expanded form, this gives

$$\begin{aligned} a_1^2 - 2(\widetilde{K}_{-1,0} + \widetilde{K}_{0,-1})a_1 &= 4\{\widetilde{K}_{00}\widetilde{K}_{-1,-1} - \widetilde{K}_{0,-1}\widetilde{K}_{-1,0} - a_0(\widetilde{K}_{00} + \widetilde{K}_{-1,-1}) + \\ &+ \widetilde{K}_{00}(\omega' - \tfrac{1}{2}\vartheta)^2 + \widetilde{K}_{-1,-1}(\omega' + \tfrac{1}{2}\vartheta)^2 + [a_0 - (\omega' - \tfrac{1}{2}\vartheta)^2][a_0 - (\omega' + \tfrac{1}{2}\vartheta)^2] \} . \end{aligned} \tag{6.2}$$

The admissible values of $\omega'$ are those that fulfill the conditions (as detailed in the previous section) for the polynomial equation (6.2), with complex coefficients, to have at least one real root in the unknown $a_1$, and those real roots are the stable values of $a_1$.

## 7. The cases of real roots of the characteristic equation.

On the sides AB and AC of the boundary of the stability region on the pq-plane (Fig. 2.1), we have solutions of the characteristic equation $(1, q)$ (on AB) and $(-1, q)$ (on AC), $q \in (-1,0) \cup (0,1)$ (the case $q = 0$ is not meaningful in the context of Floquet exponents). We set $\lambda = -\dfrac{\ln|q|}{T}$ (thus $\lambda > 0$). Then the sides AB and AC correspond to pairs of solutions of the integro-differential equation of section 3

$$x_1(t) = y_1(t),\ x_2(t) = \exp(-\lambda t) y_2(t),\ y_1, y_2\ \text{periodic},\ \lambda > 0 ;$$
$$\tilde{x}_1(t) = \tilde{y}_1(t),\ \tilde{x}_2(t) = \exp(-\lambda t)\tilde{y}_2(t),\ \tilde{y}_1, \tilde{y}_2\ \text{antiperiodic},\ \lambda > 0 .$$

In these cases, the infinite matrices $\mathbf{D}(\vartheta,\gamma)$ and $\tilde{\mathbf{D}}(\vartheta,\gamma)$ with real $\gamma$ are Hermitian, and so are appropriate finite truncations of these matrices. We shall utilize the following property of Hermitian matrices: let $A_0, A_1, ..., A_k$ be Hermitian matrices of the same dimensions, and let $x_1, x_2, ..., x_k$ be real variables; then the determinant of the Hermitian matrix $A_0 + \sum_{i=1}^{k} x_i A_i$ is a polynomial with real coefficients in the variables $x_1, x_2, ..., x_k$.

For the case of a pair of solutions, one of which is periodic and the other damped periodic, we have the system

$$\det(\mathbf{D}(\vartheta,0)) = \det(\mathbf{D}(\vartheta,\lambda)) = 0 .$$

Here, we can observe phenomena that have nothing in common with what happens in the ODE case.

The first-order approximation of the determinantal equations is the system

$$K_{00}(\vartheta,0) - a_0 = K_{00}(\vartheta,\lambda) - a_0 = 0 . \tag{7.1}$$

This system gives, in general, values of both $a_0$ and $\vartheta$. The value(s) of $\vartheta$ are determined from the equation

$$K_{00}(\vartheta,0) - K_{00}(\vartheta,\lambda) = 0 , \tag{7.2}$$

and the quasi-static value of $a_0$ is

$$a_0^* = K_{00}(\vartheta,0) . \tag{7.3}$$

The third-order approximation is the system

$$\begin{vmatrix} K_{-1,-1}(\vartheta,0)+\vartheta^2-a_0 & K_{-1,0}(\vartheta,0)-\tfrac{1}{2}a_1 & K_{-1,1}(\vartheta,0) \\ \overline{K}_{-1,0}(\vartheta,0)-\tfrac{1}{2}a_1 & K_{00}(\vartheta,0)-a_0 & K_{01}(\vartheta,0)-\tfrac{1}{2}a_1 \\ \overline{K}_{-1,1}(\vartheta,0) & \overline{K}_{01}(\vartheta,0)-\tfrac{1}{2}a_1 & K_{11}(\vartheta,0)+\vartheta^2-a_0 \end{vmatrix}=0 \; ;$$

$$\begin{vmatrix} K_{-1,-1}(\vartheta,-\lambda)+\vartheta^2-a_0 & K_{-1,0}(\vartheta,-\lambda)-\tfrac{1}{2}a_1 & K_{-1,1}(\vartheta,-\lambda) \\ \overline{K}_{-1,0}(\vartheta,-\lambda)-\tfrac{1}{2}a_1 & K_{00}(\vartheta,-\lambda)-a_0 & K_{01}(\vartheta,-\lambda)-\tfrac{1}{2}a_1 \\ \overline{K}_{-1,1}(\vartheta,-\lambda) & \overline{K}_{01}(\vartheta,-\lambda)-\tfrac{1}{2}a_1 & K_{11}(\vartheta,-\lambda)+\vartheta^2-a_0 \end{vmatrix}=0 \; .$$

(7.4)

For arbitrary but fixed $a_0$ and $\lambda$, the variable $a_1$ can be eliminated from the system (7.4): each of the two equations in (7.4) is a polynomial equation in $a_1$ and setting the resultant of the two equations in (7.4) equal to 0 gives an equation in $\vartheta$. Then, for each value of $\vartheta$ that solves the resultant equation, the value(s) of $a_1$ are determined by solving the 2 equations in (7.4).

The above considerations imply that, assuming the equation for $\vartheta$ has real solutions, say $\vartheta_k$, k = 1, 2, ... (either a finite or infinite number of solutions), for each value of $(a_0,\lambda)$ we have on the $(a_1,\vartheta)$- plane isolated points (isolated values of $a_1$ on each line $\vartheta=\vartheta_k$). For those positive values of $\lambda$ for which such real values of $a_1$ exist, the values of $a_1$ as $\lambda$ is treated as a continuous parameter, form arcs on the $(a_1,\vartheta)$- plane.

In an analogous way, for a pair of an anti-periodic solution and a damped anti-periodic solution, in the second-order truncation of the corresponding determinantal equations, and also taking into account the Hermitian properties of the coefficients $\widetilde{K}_{nm}$, we get the system

$$[a_0-\tfrac{1}{4}\vartheta^2-\widetilde{K}_{-1,-1}(\vartheta,0)][a_0-\tfrac{1}{4}\vartheta^2-\overline{\widetilde{K}}_{-1,-1}(\vartheta,0)]-[\tfrac{1}{2}a_1-\widetilde{K}_{0,-1}(\vartheta,0)][\tfrac{1}{2}a_1-\overline{\widetilde{K}}_{0,-1}(\vartheta,0)]=0 \; ;$$

$$[a_0-\tfrac{1}{4}\vartheta^2-\widetilde{K}_{-1,-1}(\vartheta,-\lambda)][a_0-\tfrac{1}{4}\vartheta^2-\overline{\widetilde{K}}_{-1,-1}(\vartheta,-\lambda)]-[\tfrac{1}{2}a_1-\widetilde{K}_{0,-1}(\vartheta,-\lambda)][\tfrac{1}{2}a_1-\overline{\widetilde{K}}_{0,-1}(\vartheta,\lambda)]$$

(7.5)

Because of the particularly simple nature of these equations, which can be written as

$$|a_0-\tfrac{1}{4}\vartheta^2-\widetilde{K}_{-1,-1}(\vartheta,0)|=|\tfrac{1}{2}a_1-\widetilde{K}_{0,-1}(\vartheta,0)| \; ;$$
$$|a_0-\tfrac{1}{4}\vartheta^2-\widetilde{K}_{-1,-1}(\vartheta,-\lambda)|=|\tfrac{1}{2}a_1-\widetilde{K}_{0,-1}(\vartheta,-\lambda)| \; ,$$

(7.6)

we can give a simple geometric interpretation to this system. It is plain that an equation $|z-z_0|=R$, on the complex plane (with $R > 0$, $z_0$ a fixed complex number, and z a variable complex number) represents a circle with center at $z_0$ and radius R. Consequently, the system (7.6) for the real-valued unknown $a_1$ expresses the requirement that the two circles

$$|z-\widetilde{K}_{0,-1}(\vartheta,0)|=|a_0-\tfrac{1}{4}\vartheta^2-\widetilde{K}_{-1,-1}(\vartheta,0)| \text{ and } |z-\widetilde{K}_{0,-1}(\vartheta,-\lambda)|=|a_0-\tfrac{1}{4}\vartheta^2-\widetilde{K}_{-1,-1}(\vartheta,-\lambda)|$$

should intersect and should have at least one of their intersection points on the real axis. This interpretation can be used as the starting point of a derivation of conditions for existence of a real root; we skip the details, as they seem to fall outside the main context of this paper.

## 8. The vertices B and C on the stability diagram.

The vertex C corresponds to periodic solutions of the integro-differential equation (3.1), and the vertex B corresponds to anti-periodic solutions.

For these points, the determinantal equations involve the determinants of Hemitian matrices. They are obtained from the equations of the previous section by setting $\lambda = 0$. Consequently, we obtain equations with real coefficients.

For the third-order approximation of the equation for periodic solutions, we have the equation

$$\begin{vmatrix} K_{-1,-1}(\vartheta,0) + \vartheta^2 - a_0 & K_{-1,0}(\vartheta,0) - \frac{1}{2}a_1 & K_{-1,1}(\vartheta,0) \\ K_{0,-1}(\vartheta,0) - \frac{1}{2}a_1 & K_{00}(\vartheta,0) - a_0 & K_{01}(\vartheta,0) - \frac{1}{2}a_1 \\ K_{1,-1}(\vartheta,0) & K_{10}(\vartheta,0) - \frac{1}{2}a_1 & K_{11}(\vartheta,0) + \vartheta^2 - a_0 \end{vmatrix} = 0 \quad . \tag{8.1}$$

By expanding this determinant, and also taking into account the Hermitian properties of the terms $K_{nm}(\vartheta,0)$, we obtain the following quadratic equation for $a_1$:

$$\begin{aligned} &(K_{00}(\vartheta,0) - a_0)\left[\left|K_{11}(\vartheta,0) + \vartheta^2 - a_0\right|^2 + 2\left|K_{-1,1}(\vartheta,0)\right|^2\right] + \\ &+ 2\Re\left[\left(K_{-1,0}(\vartheta,0) - \frac{1}{2}a_1\right)\left(K_{01}(\vartheta,0) - \frac{1}{2}a_1\right)\overline{K}_{-1,1}(\vartheta,0)\right] - \\ &- \left[\left|K_{-1,0}(\vartheta,0) - \frac{1}{2}a_1\right|^2 + \left|K_{01}(\vartheta,0) - \frac{1}{2}a_1\right|^2\right]\Re\bigl(K_{11}(\vartheta,0) + \vartheta^2 - a_0\bigr) = 0 \quad . \end{aligned} \tag{8.2}$$

One qualitative observation is that (8.2) is an equation, in $a_1$, that changes degree (from second to first) at values of $\vartheta$ at which the coefficient of $\left(\frac{1}{2}a_1\right)^2$, namely $2\Re(K_{-1,1}(\vartheta,0) - (K_{11}(\vartheta,0) + \vartheta^2 - a_0))$, vanishes; in other words, there is a sort of bifurcation phenomenon at such values of $\vartheta$. Such values of $\vartheta$ may be either isolated points or intervals. At isolated points where the degree of the polynomial equation (8.2) changes, the curves may merge so as to be continuous at that point, or they may have vertical asymptotes.

For anti-periodic solutions, we look at the second order approximation:

$$[a_0 - \tfrac{1}{4}\vartheta^2 - \widetilde{K}_{-1,-1}(\vartheta,0)][a_0 - \tfrac{1}{4}\vartheta^2 - \overline{\widetilde{K}}_{-1,-1}(\vartheta,0)] - [\tfrac{1}{2}a_1 - \widetilde{K}_{0,-1}(\vartheta,0)][\tfrac{1}{2}a_1 - \overline{\widetilde{K}}_{0,-1}(\vartheta,0)] = 0 \tag{8.3}$$

which has the solution

$$\frac{1}{2}a_1 = \Re\big(\widetilde{K}_{0,-1}(\vartheta,0)\big) \pm \sqrt{\left|a_0 - \frac{1}{4}\vartheta^2 - \widetilde{K}_{-1,-1}(\vartheta,0)\right|^2 - \big(\Im(\widetilde{K}_{0,-1}(\vartheta,0))\big)^2} \qquad (8.4)$$

provided, of course, the expression inside the square root is nonnegative.

The fourth-order truncation for anti-periodic solutions tields a polynomial equation, for the unknown $a_1$, with leading term $\left(\frac{1}{2}a_1\right)^4$, so there is no phenomenon of changing degree of the equation for ant-periodic solutions. By the same type of argument as the second justification given in section 10 below (following eq. (10.12)), the leading term in the expansion of the truncation determinanatal equation at every even order, say 2N, for anti-periodic solutions, has leading term $\left(\frac{a_1}{2}\right)^{2N}$, and thus there is never a phenomenon of changing degree for the anti-periodic solutions.

## 9. The vertex A of the stability diagram.

The vertex A corresponds to the coexistence of a periodic and an anti-periodic solution of the second-order integro-differential equation. Therefore, at this point we have a system of two determinantal equations with real coefficients for powers of $a_0$ and $a_1$. For arbitrary but fixed $a_0$, this system gives, in general, isolated values for the pair $(a_1, \vartheta)$.

At the level of third-order truncation for the periodic solutions and second-order truncation for the anti-periodic solutions, the system becomes

$$[K_{-1,-1}(\vartheta,0) + \vartheta^2 - a_0][K_{00}(\vartheta,0) - a_0][K_{11}(\vartheta,0) + \vartheta^2 - a_0] +$$
$$+ [K_{-1,0}(\vartheta,0) - \tfrac{1}{2}a_1]K_{01}(\vartheta,0)\overline{K}_{-1,1}(\vartheta,0) + [\overline{K}_{-1,0}(\vartheta,0) - \tfrac{1}{2}a_1][\overline{K}_{01}(\vartheta,0) - \tfrac{1}{2}a_1]K_{-1,1}(\vartheta,0) -$$
$$- K_{-1,1}(\vartheta,0)\overline{K}_{-1,1}(\vartheta,0)[K_{00}(\vartheta,0) - a_0] - [K_{01}(\vartheta,0) - \tfrac{1}{2}a_1][\overline{K}_{01}(\vartheta,0) - \tfrac{1}{2}a_1][K_{-1,-1}(\vartheta,0) + \vartheta^2 - a_0] +$$
$$+ [K_{-1,0}(\vartheta,0) - \tfrac{1}{2}a_1][K_{0,1}(\vartheta,0) - \tfrac{1}{2}a_1]\overline{K}_{-1,1}(\vartheta,0) = 0 \;;$$
$$[a_0 - \tfrac{1}{2}\vartheta^2 - \widetilde{K}_{-1,-1}(\vartheta,0)][a_0 - \tfrac{1}{2}\vartheta^2 - \overline{\widetilde{K}}_{-1,-1}(\vartheta,0)] - [\tfrac{1}{2}a_1 - \widetilde{K}_{0,-1}(\vartheta,0)][\tfrac{1}{2}a_1 - \overline{\widetilde{K}}_{0,-1}(\vartheta,0)] = 0 \quad .$$

(9.1)

The system of the two truncated equations is a system of two polynomial equations, with real coefficients, of the same degree in $a_1$ (for the case of third- and second-order truncations, as in (9.1), a system of two quadratic equations in $a_1$). By setting the eliminant of these two quadratic equations equal to 0, we obtain an equation in $\vartheta$. The real solutions (if any) of that equation give the discrete values of $\vartheta$ that correspond to the vertex A of the stability diagram. For each such value of $\vartheta$, the common real roots (with $a_1$ as unknown) of the two polynomial equations give the corresponding discrete values of $a_1$.

In the case considered in (9.1), things are slightly simplified because of the particular form of the second equation in (9.1), namely $|\tfrac{1}{2}a_1 - \widetilde{K}_{0,-1}(\vartheta,0)| = |a_0 - \tfrac{1}{2}\vartheta^2 - \widetilde{K}_{-1,-1}(\vartheta,0)|$, which has the real solutions (for $a_1$) given in the previous section.

## 10. Kernels of exponential type.

As mentioned previously, the absolute convergence of those integrals that correspond to the boundary of the stability domain is achieved for kernels that have exponential decay in the variable $t-s$. For this reason, we shall investigate the Strutt chart for kernels consisting of exponential terms. This class of kernels admits certain detailed calculations in exact form without any numerical approximations, and can therefore shed light on the nature of the network of curves that can appear on the Strutt chart.

We shall now consider kernels of the form

$$K(t,s) := \sum_{\alpha} c_{\alpha} \exp(-\mu_{\alpha}(t-s)) \tag{10.1}$$

as well as the particular case

$$K(t,s) := c \exp(-\mu(t-s)) \ . \tag{10.2}$$

The idea for kernels of the form (10.1) comes from [BA], a paper that deals with the modelling of viscoelastic materials ( but not with problems of dynamic stability).

In (10.1), we postulate

$$c_{\alpha} > 0, \ \mu_{\alpha} > 0, \ \sum_{\alpha} c_{\alpha} < \infty . \tag{10.3}$$

Then

$$K_{nm}(\vartheta,\gamma) = \sum_{\alpha} \frac{c_{\alpha}}{\mu_{\alpha} + \gamma + im\vartheta} \delta_{nm} . \tag{10.4}$$

With $\gamma = 0$, we have

$$K_{nm}(\vartheta,0) = \sum_{\alpha} \frac{c_{\alpha}(\mu_{\alpha} - im\vartheta)}{\mu_{\alpha}^{2} + (m\vartheta)^{2}} \delta_{nm} \ ; \ K_{-n,-m}(\vartheta,0) = \overline{K}_{nm}(\vartheta,0) = \sum_{\alpha} \frac{c_{\alpha}(\mu_{\alpha} + mi\vartheta)}{\mu_{\alpha}^{2} + (m\vartheta)^{2}} \delta_{nm} \tag{10.5}$$

and consequently

$$K_{nm}(\vartheta,0) + K_{-n,-m}(\vartheta,0) = 2\sum_{\alpha} \frac{c_{\alpha}\mu_{\alpha}}{\mu_{\alpha}^2 + (m\vartheta)^2}\delta_{nm} ;$$

$$K_{nm}(\vartheta,0) K_{n'm'}(\vartheta,0) = \sum_{\alpha,\alpha'} \frac{c_{\alpha}c_{\alpha'}[\mu_{\alpha}\mu_{\alpha'} + mm'\vartheta^2 + i(m'\mu_{\alpha} - m\mu_{\alpha'})\vartheta]}{(\mu_{\alpha}^2 + (m\vartheta)^2)(\mu_{\alpha'}^2 + (m'\vartheta)^2)} .$$

(10.6)

On the basis of these results, we can find the equations for the third-order approximation to the curves (on the $(\vartheta, a_1)$ plane) that correspond to periodic solutions of the integro-differential equation of section 3. For simplicity, we write down the equations for the case of the simpler model (10.2).

After some calculations, the third-order approximation becomes

$$\left(\frac{c}{\mu} - a_0\right)\left(\frac{c^2 + 2c\mu(\vartheta^2 - a_0)}{\mu^2 + \vartheta^2} + (\vartheta^2 - a_0)^2\right) - \frac{a_1^2}{2}\left(\frac{c\mu}{\mu^2 + \vartheta^2} + \vartheta^2 - a_0\right) = 0 \qquad (10.7)$$

which gives

$$a_1^2 = 2\left(\frac{c}{\mu} - a_0\right)\left[\vartheta^2 - a_0 + \frac{c^2 + c\mu(\vartheta^2 - a_0)}{c\mu + (\vartheta^2 - a_0)(\mu^2 + \vartheta^2)}\right] \qquad (10.8)$$

This equation is meaningful for values of the parameters that make the right-hand side nonnegative; thus the admissible values of $\vartheta^2$ in (10.8) depend on $a_0$, $c$, and $\mu$.

The curve on the $(\vartheta^2, a_1^2)$- plane represented by (10.8) has a vertical asymptote near positive roots (in the variable $\vartheta^2$) of the denominator $c\mu + (\vartheta^2 - a_0)(\mu^2 + \vartheta^2)$; such positive roots exist if $a_0$ is sufficiently large, and then, of course, the corresponding values of $\vartheta^2$ are

$$\vartheta_*^2 = \frac{1}{2}\left\{a_0 - \mu^2 \pm \sqrt{(a_0 - \mu^2)^2 - (2c\mu)^2}\right\} . \qquad (10.9)$$

This creates another limitation on the admissible values of $\vartheta^2$, since for $\vartheta^2$ near those roots the coefficient $a_1^2$ as given by (10.8) becomes arbitrarily large; of course, all curves obtained by the method of determinants are physically meaningful for "small" values of the relevant parameters.

Next, we examine the second-order approximation to the equations that correspond to anti-periodic solutions. The terms $\widetilde{K}_{nm}(\vartheta,\gamma)$ are calculated as

$$\widetilde{K}_{nm}(\vartheta,\gamma) = \sum_{\alpha} \frac{c_{\alpha}}{\mu_{\alpha} + \gamma + \frac{i(2m+1)}{2}\vartheta} \delta_{nm} \ . \tag{10.10}$$

Based on (10.10), we can find the second- and fourth-order approximations to the determinantal equation for anti-periodic solutions, for kernels of the simpler form (10.2).
After some calculations, the second-order approximation gives

$$a_1^2 = \frac{4}{\left(\mu^2 + \frac{\vartheta^2}{4}\right)^2} \left[ \left\{ \left(a_0 - \frac{\vartheta^2}{4}\right)\left(\mu^2 + \frac{\vartheta^2}{4}\right) - c\mu \right\}^2 + \frac{c^2 \vartheta^2}{4} \right] . \tag{10.11}$$

The right-hand side of (10.11), as a function of $\vartheta^2$, is a rational function whose denominator never vanishes, and so we have no asymptotes.

For the fourth-order approximation, we find the equation

$$\frac{1}{16}a_1^4 - \frac{1}{4}a_1^2 \left[ \left( a_0 - \frac{9\vartheta^2}{4} - \frac{c\mu}{\mu^2 + \frac{9\vartheta^2}{4}} \right)^2 + \frac{9c\frac{\vartheta^2}{4}}{\left(\mu^2 + \frac{9\vartheta^2}{4}\right)^2} + \right.$$

$$+ 2\left(a_0 - \frac{9\vartheta^2}{4} - \frac{c\mu}{\mu^2 + \frac{9\vartheta^2}{4}}\right)\left(a_0 - \frac{\vartheta^2}{4} - \frac{c\mu}{\mu^2 + \frac{\vartheta^2}{4}}\right) + \left(a_0 - \frac{9\vartheta^2}{4} - \frac{c\mu}{\mu^2 + \frac{9\vartheta^2}{4}}\right)\frac{c\vartheta}{\mu^2 + \frac{\vartheta^2}{4}} +$$

$$+ \left(a_0 - \frac{\vartheta^2}{4} - \frac{c\mu}{\mu^2 + \frac{\vartheta^2}{4}}\right)\frac{3c\vartheta}{\mu^2 + \frac{9\vartheta^2}{4}} \right] +$$

$$+ \left[ \left( \left(a_0 - \frac{9\vartheta^2}{4} - \frac{c\mu}{\mu^2 + \frac{9\vartheta^2}{4}}\right)^2 + \left(\frac{\frac{3c\vartheta}{2}}{\mu^2 + \frac{9\vartheta^2}{4}}\right)^2 \right) \left( \left(a_0 - \frac{\vartheta^2}{4} - \frac{c\mu}{\mu^2 + \frac{\vartheta^2}{4}}\right)^2 + \left(\frac{\frac{c\vartheta}{2}}{\mu^2 + \frac{\vartheta^2}{4}}\right)^2 \right) \right] = 0 \ . \tag{10.12}$$

We observe that whenever (10.12) has real roots (in the unknown $a_1^2$), they are both positive (since the zero-degree term on the left-hand side of (10.12) is always positive), and we get two curves on the $(\vartheta^2, a_1^2)$ plane. Asymptotes do not exist since all denominators which can arise are powers with rational exponents of the positive terms $\mu^2 + \frac{\vartheta^2}{4}$ and $\mu^2 + \frac{9\vartheta^2}{4}$.

The property of lack of singularities, in the anti-periodic case with exponential kernels, also holds for all even-order approximations to the determinantal equation, under quite general conditions. This can be seen as follows: in the approximation of order 2N, we have the determinant of a tridiagonal matrix in which the diagonal terms are expressions of the form

$$a_0 - \left(\frac{(2n+1)\vartheta}{2}\right)^2 - \frac{\mu \pm i\left(\frac{(2n+1)\vartheta}{2}\right)}{\mu^2 + \left(\frac{(2n+1)\vartheta}{2}\right)^2},$$

whereas the off-diagonal terms are equal to $\frac{a_1}{2}$; in this case, the expansion of the determinant has a leading term in $a_1$ whose absolute value is $\left(\frac{a_1}{2}\right)^{2N}$.

This can be shown either by induction using the general theory of recursive methods for the evaluation of a tridiagonal determinant (such as, e.g., the method of continuants [M]), or directly, by observing that, in the tridiagonal matrix under consideration, we can find the term $\frac{a_1}{2}$ in 2N locations, such that these locations have no row-index in common and no column-index in common, and this is the only set of such locations whose cardinality is 2N and each location contains $a_1$. Consequently, the determinantal equation becomes a polynomial equation with leading-term coefficient equal to a non-zero constant, and remaining coefficients that have no singularities as functions of $\vartheta$. Under these conditions, the theorem on analytic dependence of polynomial roots on the coefficients of a polynomial applies. (See, for example, [Br]. Note that the results of [Br] apply to a monic polynomial, although this is not explicitly stated in the individual enunciations of the theorems in that paper, but it follows from the overall context of the paper.) Consequently: the distinct roots (including the possibility of distinct roots some of which are multiple roots) of the resulting polynomial equation are analytic functions of the coefficients, as long as the coefficients remain in a domain such that the multiplicity of each those distinct roots remains constant. Under these stipulations, the resulting solutions for $a_1$ have no singularities in $\vartheta$, since the coefficients themselves have no singularities in $\vartheta$.

Appendix I. Existence and uniqueness of solutions of an integro-differential system of the type considered in this paper.

Although there is a substantial amount of literature on integral equations over infinite intervals (cf. [AOR]), the specific results we need for the type of integro-differential equations considered in this paper do not exist in the existing research literature. For this reason, we proved the necessary existence and uniqueness results directly from first principles.

We consider the integro-differential equation

$$\dot{x}(t) + A(t)x(t) = \int_{s=-\infty}^{t} K(t,s)x(s)\,ds\,;\quad x(0) = x_0 \tag{I.1}$$

The unknown function $x(\cdot)$ takes values in a finite-dimensional Euclidean space, and the matrix-valued functions A and K have dimensions compatible with the dimension of x.

In order to show existence and uniqueness of solutions for $-\infty < t < +\infty$, it is plain that the crucial issue is to show existence and uniqueness of solutions over $-\infty < t \leq 0$: indeed, if $x_1(t)$ is a solution of (I.1) for $-\infty < t \leq 0$, then, for t>0, we have to solve the problem

$$\dot{x}(t) + A(t)x(t) = \int_{s=-\infty}^{0} K(t,s)x_1(s)\,ds + \int_{s=0}^{t} K(t,s)x(s)\,;\quad x(0) = x_0 \tag{I.2}$$

Eq. (I.2) is a Volterra integro-differential equation, for which existence and uniqueness of solutions can be proved by standard methods. On the other hand, for t<0, eq. (I.1) is not reducible to a Volterra integral equation, but rather to a hybrid type of integral equation that contains both Volterra and Fredholm integral operators. The integral formulation of (I.1), for t<0, is

$$x(t) = x_0 + \int_{s=t}^{0} A(s)x(s)\,ds - \int_{s_1=t}^{0} \int_{s=-\infty}^{s_1} K(s_1,s)x(s)\,ds\,ds_1$$

and after a change in the order of integration we get

$$x(t) = x_0 + \int_{s=t}^{0} A(s)x(s)\,ds - \int_{s=-\infty}^{0} K_1(t,s)x(s)\,ds\,;\quad K_1(t,s) := \int_{s_1=\max(t,s)}^{0} K(s_1,s)\,ds_1 \tag{I.3}$$

Our proof of existence and uniqueness of solutions of (I.3) uses arguments and techniques that are not specific to linear equations, but rather can be applied to a more general nonlinear problem of the form

$$x(t) = z(t) + \int_{s=t}^{0} f(t,s,x(s))ds - \int_{s=-\infty}^{0} \int_{s_1=\max(t,s)}^{0} g(t,s_1,s,x(s))ds_1\, ds$$

(I.4)

We note that, when the function $z(t)$ is a constant $x_0$ and the functions $f$ and $g$ are independent of t, say $f(t,s,x) \equiv f_1(s,x)$, $g(t,s_1,s,x) \equiv g_1(s_1,s,x)$, then (I.4) is the integral formulation, for t<0, of the integro-differential initial value problem

$$\dot{x}(t) + f_1(t,x(t)) = \int_{s=-\infty}^{t} g_1(t,s,x(s))ds\,;\; x(0) = x_0$$

(I.5)

Therefore, we present the proof for the more general problem (I.4), and then explain how it can be specialized to the problem (I.3).

We assume:

<u>(A1).</u> The functions f and g are continuous in all their arguments and satisfy

$|f(t,s,x) - f(t,s,y)| \leq C_f\,|x-y|$, for all $-\infty < s \leq t \leq 0$, x, y in $\mathbb{R}^n$ ;

$|g(t,s_1,s,x) - g(t,s_1,s,y)| \leq C_g \exp(-\mu(s_1-s))|x-y|$, for all $t \leq 0$, $s \leq 0$, $\max(t,s) \leq s_1 \leq 0$, x, y in $\mathbb{R}^n$ .

<u>(A2).</u> For every $t \in (-\infty,0]$, the function $|g(t,s_1,s,0)|$ is integrable with respect to $(s_1,s)$ over the domain $\{(s_1,s): \max(t,s) \leq s_1 \leq 0,\; -\infty < s \leq 0\}$.

<u>(A3).</u> The function $z(\cdot)$ is continuous and bounded on $(-\infty,0]$.

<u>(A4).</u> The constants $\mu$, $C_f$, $C_g$ satisfy $0 < \mu - C_f < \sqrt{2}\,C_g$ .

For every $\gamma > 0$, the class $C_\gamma(-\infty, 0)$ consists of functions $w(\cdot)$ that are continuous on $(-\infty, 0]$ and have the property that $e^{\gamma t} w(t)$ is bounded on $(-\infty, 0]$.

We shall prove:

Theorem I.1. Under assumptions (A1) through (A4), there is a $\gamma^* > 0$ such that (I.4) has a unique solution in $C_{\gamma^*}(-\infty, 0)$.

Proof: For an as yet unspecified $\gamma$ that satisfies $0 < \gamma < \mu$, we set $y(t) = e^{\gamma t} x(t)$. Then y satisfies

$$y(t) = z(t)e^{\gamma t} + \int_{s=t}^{0} f(t, s, e^{-\gamma s} y(s)) \, ds - \int_{s=-\infty}^{0} \int_{s_1=\max(t,s)}^{0} g(t, s_1, s, e^{-\gamma s} y(s)) \, ds_1 \, ds \, .$$

(I.6)

For every continuous bounded function $z(\cdot)$, we set

$$(Sz)(t) := e^{\gamma t} \left[ z(t) + \int_{s=t}^{0} f(t, s, e^{-\gamma s} z(s)) \, ds - \int_{s=-\infty}^{0} \int_{s_1=\max(t,s)}^{0} g(t, s_1, s, e^{-\gamma s} z(s)) \, ds_1 \, ds \right].$$

The double integral over an infinite domain, in the equation above, is well-defined by assumptions (A1) and (A2).

For any two continuous bounded functions $y_1(\cdot)$, $y_2(\cdot)$ on $(-\infty, 0]$, we have

$$|(Sy_1)(t) - (Sy_2)(t)| \leq C_f \int_{s=t}^{0} \exp(\gamma(t-s)) \, ds [\max_{t \leq s \leq 0} |y_1(s) - y_2(s)|] +$$

$$+ C_g \int_{s=-\infty}^{0} \int_{s_1=\max(t,s)}^{0} \exp(\gamma(t-s) - \mu(s_1-s)) \, ds_1 \, ds \, [\sup_{-\infty \leq s \leq 0} |y_1(s) - y_2(s)|] \leq$$

$$\leq \frac{C_f}{\gamma} [\max_{t \leq s \leq 0} |y_1(s) - y_2(s)|] + \frac{C_g}{\gamma(\mu-\gamma)} [\sup_{-\infty \leq s \leq 0} |y_1(s) - y_2(s)|] \, .$$

(I.7)

The last inequality, in the string of inequalities above, results from explicit calculation of the integrals contained in (I.7); we omit the details.

Consequently, S will be a contraction on the space of continuous bounded functions on $(-\infty,0]$ with the supremum norm if there exits a $\gamma^*$ in $(0,\mu)$ such that

$$\frac{C_f}{\gamma^*} + \frac{C_g}{\gamma^*(\mu-\gamma^*)} < 1 \ .$$

It is a matter of elementary algebra to verify that condition (A4) suffices for the existence of such a $\gamma^*$. By the basic fixed point theorem for contraction mappings, it follows that eq. (I.6) has a unique solution in the space of continuous bounded functions on $(-\infty,0]$ with the supremum norm, and consequently eq. (I.4) has a unique solution in $C_{\gamma^*}(-\infty,0)$ . ///

The integral formulation of the integro-differential problem (I.5) for t>0 is

$$x(t) = x_0 - \int_{s=0}^{t} f_1(s,x(s))\,ds + \int_{s=-\infty}^{t} \int_{s_1=\max(s,0)}^{t} g_1(s_1,s,x(s))\,ds_1\,ds \tag{I.8}$$

Therefore, we consider the more general problem

$$x(t) = z(t) - \int_{s=0}^{t} f(t,s,x(s))\,ds + \int_{s=-\infty}^{t} \int_{s_1=\max(s,0)}^{t} g(t,s_1,s,x(s))\,ds_1\,ds \tag{I.9}$$

where $x(s)$, $s \leq 0$, is the unique solution of (I.4) as obtained in theorem 1.

We shall postulate the following property:

(A5). The conditions (A1) through (A3) also hold for all $(t,s_1,s)$ that satisfy $t \geq 0$, $-\infty < s \leq t$, $\max(s.0) \leq s_1 \leq t$.

Then we have:

<u>Theorem I.2.</u> We assume conditions (A1) through (A5). Let $x_1(t)$, $t<0$, be the solution obtained in theorem 1, and $\gamma^*$ the value of $\gamma$ specified in theorem 1. Then the solution $x_1$ can be extended to a solution $x(t)$ of (I.9) for $t>0$.

<u>Proof:</u> We write (I.9) in the equivalent form

$$x(t) = z(t) + \int_{s=-\infty}^{0} \int_{s_1=0}^{t} g(t,s_1,s,x_1(s))ds_1\,ds - \int_{s=0}^{t} f(t,s,x(s))ds +$$

$$+ \int_{s=0}^{t} g_1(t,s,x(s))ds \; ; \; g_1(t,s,x) := \int_{s_1=s}^{t} g(t,s_1,s,x)ds_1$$

(I.10)

By assumptions (A1) and (A2) and the fact that $0 < \gamma^* < \mu$, the integral
$\int_{s=-\infty}^{0} \int_{s_1=0}^{t} g(t,s_1,s,x_1(s))ds_1\,ds$ is absolutely convergent and defines a continuous function of t.

The function $g_1$ is estimated as follows, for $s>0$:

$$|g_1(t,s,w_1) - g(t,s,w_2)| \leq C_g \int_{s_1=s}^{t} \exp(-\mu(s_1-s))|w_1-w_2|ds =$$

$$= \frac{C_g}{\mu}[1-\exp(-\mu(t-s))]|w_1-w_2| \leq \frac{C_g}{\mu}|w_1-w_2|, \; \forall \; w_1,w_2 \text{ in } \mathbb{R}^n$$

(I.11)

Consequently, (I.10) is a Volterra integral equation with continuous forcing term
$z(t) + \int_{s=-\infty}^{0} \int_{s_1=0}^{t} g(t,s_1,s,x_1(s))ds_1\,ds$ and Lipschitz (nonlinear) kernel $g_1 - f$, and existence and uniqueness of continuous solutions for $t>0$ follows from the standard theory of Volterra integral equations. ///